\documentclass[11pt]{article}
\usepackage[ansinew]{inputenc}
\usepackage{graphicx}
\usepackage{color}

\usepackage{enumerate,latexsym}
\usepackage{latexsym}
\usepackage{amsmath,amssymb}
\usepackage{graphicx}
\usepackage{times}

\newfont{\bb}{msbm10 at 11pt}
\newfont{\bbsmall}{msbm8 at 8pt}
\def\rth{\mathbb{R}^3}
\def\R{\mathbb{R}}

\def\Z{\mathbb{Z}}
\def\C{\mathbb{C}}

\def\D{\mathbb{D}}

\def\esf{\mathbb{S}}

\newcommand{\ben}{\begin{enumerate}}
\newcommand{\bit}{\begin{itemize}}
\newcommand{\een}{\end{enumerate}}
\newcommand{\eit}{\end{itemize}}

\newcommand{\Int}{\mbox{Int}}

\newcommand{\Sov}{\overline{\Sigma}}
\newcommand{\Dhat}{\widehat{D}}

\def\a{{\alpha}}

\def\g{{\gamma}}

\def\de{{\delta}}
\def\be{{\beta}}
\def\ve{{\varepsilon}}

\def\centerbmp#1#2#3{\vskip#2\relax\centerline{\hbox to#1{\special
    {bmp:#3 x=#1, y=#2}\hfil}}}

\newtheorem{theorem}{Theorem}[section]

\newtheorem{remark}[theorem]{Remark}
\newtheorem{corollary}[theorem]{Corollary}

\newtheorem{assertion}[theorem]{Assertion}

\newenvironment{proof}{\smallskip\noindent{\it Proof.}\hskip \labelsep}
{\hfill\penalty10000\raisebox{-.09em}{$\Box$}\par\medskip}

\begin{document}
\begin{title}
{Finite type annular ends for harmonic functions}
\end{title}
\vskip .2in

\begin{author}
{William H. Meeks III\thanks{This material is based upon
   work for the NSF under Award No. DMS-1309236.
   Any opinions, findings, and conclusions or recommendations
   expressed in this publication are those of the authors and do not
   necessarily reflect the views of the NSF.}
   \and Joaqu\'\i n P\' erez\thanks{Research supported in part
by the MINECO/FEDER grant no. MTM2014-52368-P.
}}
\end{author}
\maketitle
\begin{abstract}
In this paper we describe the notion of an annular end of a Riemann
surface being of finite type with respect to some harmonic function
and  prove some theoretical results relating the conformal structure
of such an annular end to the level sets of the harmonic function.
As an application of these results, we obtain important information
on the conformal type of any properly immersed minimal surface $M$
in $\rth$ with compact boundary and which intersects some plane
transversely in a finite number of arcs; in particular, such an $M$
is a parabolic Riemann surface. This information is applied by the
authors in~\cite{mpe3} to classify the asymptotic behavior of
annular ends of a complete embedded minimal annulus with compact
boundary in terms of the flux vector along its boundary. In the
present paper, we apply this information to understand and
characterize properly immersed minimal surfaces in $\mathbb{R}^3$ of
finite total curvature, in terms of their intersections with two
nonparallel planes.

\vspace{.3cm}

\noindent{\it Mathematics Subject Classification:} Primary 53A10,
   Secondary 49Q05, 53C42

\noindent{\it Key words and phrases:} Minimal surface, finite type,
harmonic function, parabolic Riemann surface, hyperbolic Riemann
surface, angular limits.
\end{abstract}

\section{Introduction.}
Given a nonconstant harmonic function $f\colon M\to \R$ on a
noncompact Riemann surface with compact boundary, we say that an
annular end\footnote{A proper subdomain $E\subset M$ is an {\it
annular end} if it is homeomorphic to $\esf^1\times [0,1)$.}
$E$ of
 $M$ has {\it finite type for} $f$ if for some $t_0\in \R$, the
one-complex $f^{-1}(t_0)\cap E$ has a finite number of ends.
Observe that if $E'$ is an annular subend of $E$, then $E'$ has
finite type for $f$ if and only if $E$ has finite type for $f$. Thus, 
in the sequel we will assume that $f$ has no critical points in the
boundary of $E$.

Note that $f^{-1}(t_0)$ might fail to be smooth: at each critical
point $p$ of $f$ lying in $f^{-1}(t_0)$, this level set consists
locally of an equiangular system of curves crossing at $p$; also
note that $f^{-1}(t_0)$ cannot bound a compact subdomain in
$E-\partial E$ by the maximum principle. These observations imply
that $f^{-1}(t_0)\cap E$ has a finite number of ends if and only if
$f^{-1}(t_0)\cap E$ has a finite number of components and a finite
number of crossing points.

\begin{figure}
\begin{center}
\includegraphics[width=6cm]{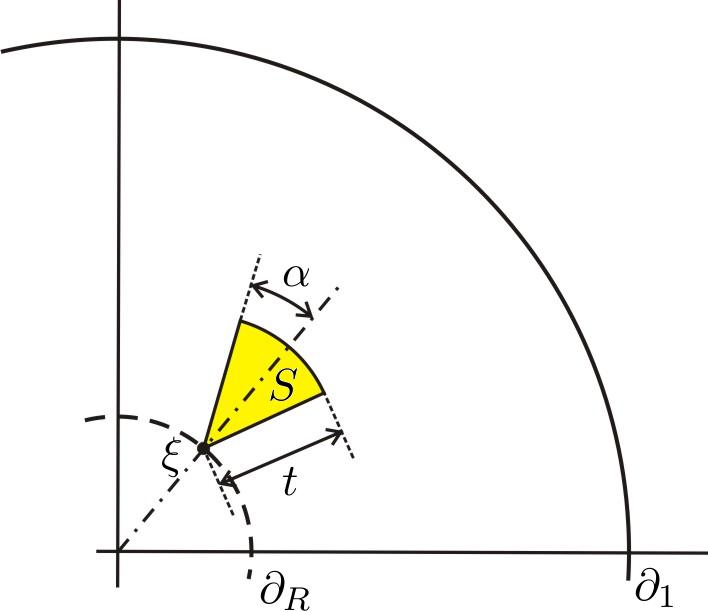}
\caption{Angular sector in $A(R,1)$, centered at
$\xi \in \partial _R$.} \label{fig1}
\end{center}
\end{figure}

When $E$ is an annular end of finite type for $f$, we will prove
several results on the level sets of $f|_E$ depending on the
conformal type of $E$. In order to describe these results, we first
fix some notation. For $R\in [0,1)$, let $A(R,1)=\{z\in \C \mid R<
|z|\leq 1\}$, $\overline{A}(R,1) =\{z\in \C \mid R\leq |z|\leq 1\}$,
$\partial_R=\{t\in\C \mid |z|=R\}$ and $\partial_1=\{ z\in\C \mid
|z|=1\}$. Thus, the closure $\overline{A}(0,1)$ in $\C$ is the closed unit
disk $\overline{\D }$ and $A(0,1)=\overline{\D}-\{ 0\} $. If $0<R<1$
and $F\colon A(R,1)\to \C $ is a holomorphic function, then $F$ is
said to have {\it angular limit $F(\xi )\in \C \cup \{ \infty \} $
at $\xi \in \partial _R$} if $\lim _{z\to \xi ,z\in S}F(z)$ exists
and equals $F(\xi )$ for every angular sector $S\subset A(R,1)$
centered at $\xi $, whose median line is the radial arc $[1,2] \xi$,
with small radius $t\in (0,1-R)$ and total angle $2\a $, $0<\a
<\frac{\pi }{2}$, see Figure~\ref{fig1}. This definition of having
angular limit can be directly extended to (real valued) harmonic
functions.

We can now state our main result.

\begin{theorem}
\label{th1.1}
Suppose $f\colon M\to \R$ is a nonconstant harmonic function
and $E$ is an annular end of $M$ of finite type for $f$. Then:
\ben
\item If $E$ is conformally $A(0,1)$, then the holomorphic one-form
$\omega=df+idf^*$ on $A(0,1)$ extends to a meromorphic one-form
$\widetilde{\omega }$ on the closed unit disk $\overline{\D}$ (here
$f^*$ denotes the locally defined conjugate harmonic function of
$f$). Furthermore:
\begin{enumerate}
\item If for some $t_0\in \R $ the one-complex $f^{-1}(t_0)\cap E$ has
$2k$ ends (note that the number of ends is always even since $f$
alternates values greater and smaller than $t_0$ in consecutive
components of the complement of $f^{-1}(t_0)\cap E$ around $z=0$) and
$\widetilde{\omega }$ has a pole at $z=0$, then this pole is of
order $k+1$. In particular, for every $t\in \R$, the level set
$f^{-1}(t)$ has exactly $2k$ ends.
\item $\widetilde{\omega }$ is holomorphic if and only if $f$
is bounded on $E$ (in which case $f$ admits a well-defined harmonic
conjugate function on $\overline{\D }$).
\item If $\int _{\a }|df|<\infty $ for some end representative
$\a $ of $f^{-1}(t_0)\cap E$, then $f$ is bounded on $E$.
\end{enumerate}
\item If $E$ is conformally $A(R,1)$ for some $R\in (0,1)$, then:
\ben
\item $f$ has angular limits almost everywhere on $\partial_R$.
\item Given $t_0 \in \R$ such that $f^{-1}(t_0)\cap E$ has a finite number of
ends, then the limit set of each end of $f^{-1}(t_0)\cap E$ consists of a
single point in $\partial_R$.  In particular, the closure in
$A(R,1)$ of at least one component of $\{z\in A(R,1) \mid f(z)\neq
t_0\}$ is hyperbolic\footnote{A noncompact Riemann surface $\Sigma$ with
boundary is {\it hyperbolic } if its boundary fails to have full
harmonic measure (equivalently, bounded harmonic functions on
$\Sigma $ are not determined by their boundary values); otherwise,
$\Sigma$ is called {\it parabolic}.}. \een \een
\end{theorem}

In the special case that the flux $\int_{\partial _1}\frac{\partial
f}{\partial r}\, ds$ vanishes, then item~2 of
Theorem~\ref{th1.1} follows from the proof of Theorem 7.1
in~\cite{mr7}.

Theorem~\ref{th1.1} is motivated by applications of it to minimal
surface theory.
%
For specific applications in~\cite{mpe3}, we will need the following
corollary to Theorem~\ref{th1.1}. We remark that
L\'opez~\cite{lop4} obtained related results of Picard type for
minimal surfaces under the assumption that the second sentence of
Corollary~\ref{c1.2} holds, and that Alarc\'on and L\'opez~\cite{allo1} also applied some of the main results in this paper. 

\begin{corollary}
\label{c1.2} Suppose $X=(x_1,x_2,x_3)\colon A(R,1)\to \rth$ is a
proper, conformal minimal immersion and that for some horizontal plane
$P\subset \R^3$, the one-complex $X^{-1}(P)$ has a finite number
$2k$ of ends. Then $R=0$ and the holomorphic height differential
$dx_3 + idx_3^*$ extends meromorphically to $\overline{\D}$ with a
pole of order $k+1$. In particular, for every horizontal plane
$P'\subset \R^3$, the level set $X^{-1}(P')$ has exactly $2k$ ends.

Furthermore, after replacing $A(0,1)$ by $A(0,R')=\{z\in \C \mid
0<|z|\leq R'\}$ for some $R' \in (0,1)$, the Gauss map of the
induced minimal immersion $X|_{A(0,R')}$ is never vertical. Hence,
on $A(0,R')$, the meromorphic Gauss map of $X$ can be expressed as
$g(z)=z^ne^{H(z)}$ for some $n\in \Z$ and for some holomorphic
function $H\colon A(0, R')\to \C$. Furthermore, the following three
statements are equivalent:
\begin{enumerate}
\item $X$ has finite total curvature.
\item $H$ is bounded on $A(0,R')$.
\item There are two nonparallel planes $P_1, P_2$ such that
each of these planes individually intersects the surface
transversely in a finite number of immersed curves (if $X$ is not
flat, then we can change ``there are two'' by ``for all'' in this
statement).
\end{enumerate}
Additionally, if the flux $\int _{\a }\frac{\partial x_3}{\partial
\eta }ds$ is finite for some end representative $\a $ of
$X^{-1}(P)$, then $X$ has finite total curvature and $X(A(R,1))$ is
asymptotic to $P$ with finite multiplicity.
\end{corollary}

\section{Preliminaries on complex analysis.}
In this section we recall the statements of two classical theorems
in the theory of functions of one complex variable that we will
apply to prove Theorem~\ref{th1.1}.

\begin{theorem}[Plessner~\cite{ple1}]
\label{thmplessner} If $F$ is a meromorphic function in the open
unit disk $\D $, then, for almost all $\xi \in \partial \D$, either
$F$ has a finite angular limit at $\xi$ or $F(S)$ is dense in $\C \cup
\{ \infty \} $, for every angular sector $S=\{ z\in \D \mid |\arg
(1-\overline{z}\xi )|<\a , |z-\xi |<t\} $ centered at $\xi $ of
radius $t$, aperture angle $2\a $ (here $0<\a <\frac{\pi}{2 }$) and
median line $[0,\xi ]$.
\end{theorem}

\begin{theorem}[Privalov~\cite{priv1}]
\label{thmPrivalov} Let $F$ be a meromorphic function in $\D $. If
$F$ has angular limit~$0$ in a subset of positive measure in
$\partial \D $, then $F$ vanishes identically.
\end{theorem}

\section{The proof of Theorem~\ref{th1.1}.}
Suppose $t_0\in \R$ and $f\colon A(R,1)\to \R$ is a nonconstant
harmonic function with $f^{-1}(t_0)$ having $n+1$ ends. The
ends of $f^{-1}(t_0)$ can be represented
by a family $\{\a_0,\a_1,\ldots, \a_n\}$ of pairwise
disjoint proper arcs in $A(R,1)$. Our first observation in proving
Theorem~\ref{th1.1} is  given in the following assertion.

\begin{assertion}
 \label{a2.1}
There exists a simple closed analytic curve $\beta\colon \esf^1\to
A(R,1)$ which is topologically parallel to $\partial_1$ and which
intersects $f^{-1}(t_0)$ transversely in  $n+1$ points $p_0=\a_0\cap
\beta, p_1=\a_1\cap \beta, \ldots, p_n=\a_n\cap \beta $.
Furthermore, after replacing  $A(R,1)$ by the closure of the
subdomain of $A(R,1)-\beta$ which is disjoint from
$\partial_1$, then $f^{-1}(t_0)$ consists of $n+1$ disjoint proper arcs
representing  the ends of $f^{-1}(t_0)$.
\end{assertion}
\begin{proof}
This assertion follows immediately from elementary separation
properties of curves and the conformal classification of annular
domains.
\end{proof}

By Assertion~\ref{a2.1}, after replacing $A(R,1)$ by a subend
and parameterizing this subend by some $A(R',1)$, we may assume that
$f\colon A(R,1)\to \R$ is harmonic and analytic (up to and including
the boundary $\partial _1$) and $f^{-1}(t_0)$ is a finite collection
$\{\a_0, \a_1, \ldots, \a_n\}$ of pairwise disjoint, properly
embedded arcs transverse to $\partial_1$, and each arc $\a_i$ has
its starting point $p_i$ in $\partial_1$, $i=0,1, \ldots, n$. We can
also assume that $\{\a_0,\a_1,\ldots, \a_n\}$ are cyclically ordered
in a counterclockwise manner.

\begin{assertion}
\label{a2.2} Suppose that $R>0$ and that the limit set
$L(f^{-1}(t_0))\subset \partial_R$ of $f^{-1}(t_0)$ is not equal to
$\partial_R$. Then, item~2 in the statement of Theorem~\ref{th1.1}
holds.
\end{assertion}
 \begin{proof}
Let $\sigma$ be a compact embedded arc in $\overline{A}(R,1)-
\bigcup^n_{i=0}\a_i$ with one end point in $\partial_1$ and the
other end point in $\partial_R$ (note that $\sigma $ exists since
$\partial _R-L(f^{-1}(t_0))\neq \mbox{\O }$). Let $\sigma(\ve)$ be a
small, open regular neighborhood of $\sigma$ in
$\overline{A}(R,1)-\bigcup^n_{i=0}\a_i$ which is at a positive
distance from $\bigcup_{i=0}^n\a_i$ and chosen so that $\partial
[A(R,1)-\sigma(\ve)]$ is a smooth connected arc. By the conformal
classification of annular domains and boundary regularity of
holomorphic functions, there exists a conformal diffeomorphism (see
Figure \ref{fig2})
\[
\eta\colon \Delta=\left\{ z=x+iy\in \C \mid |z|<1,\ y \geq
0\right\}\to [A(R,1)-\sigma(\ve)].
\]

Consider the induced harmonic function $\widehat{f}=f\circ \eta
\colon \Delta\to \R$. Since $\Delta$ is simply connected,
$\widehat{f}$ admits a well-defined harmonic conjugate function
$\widehat{f}^*$. Hence, the function
$F=\widehat{f}+i\widehat{f}^*\colon \Delta\to\C$ is holomorphic.
\begin{figure}
\begin{center}
\includegraphics[height=5cm]{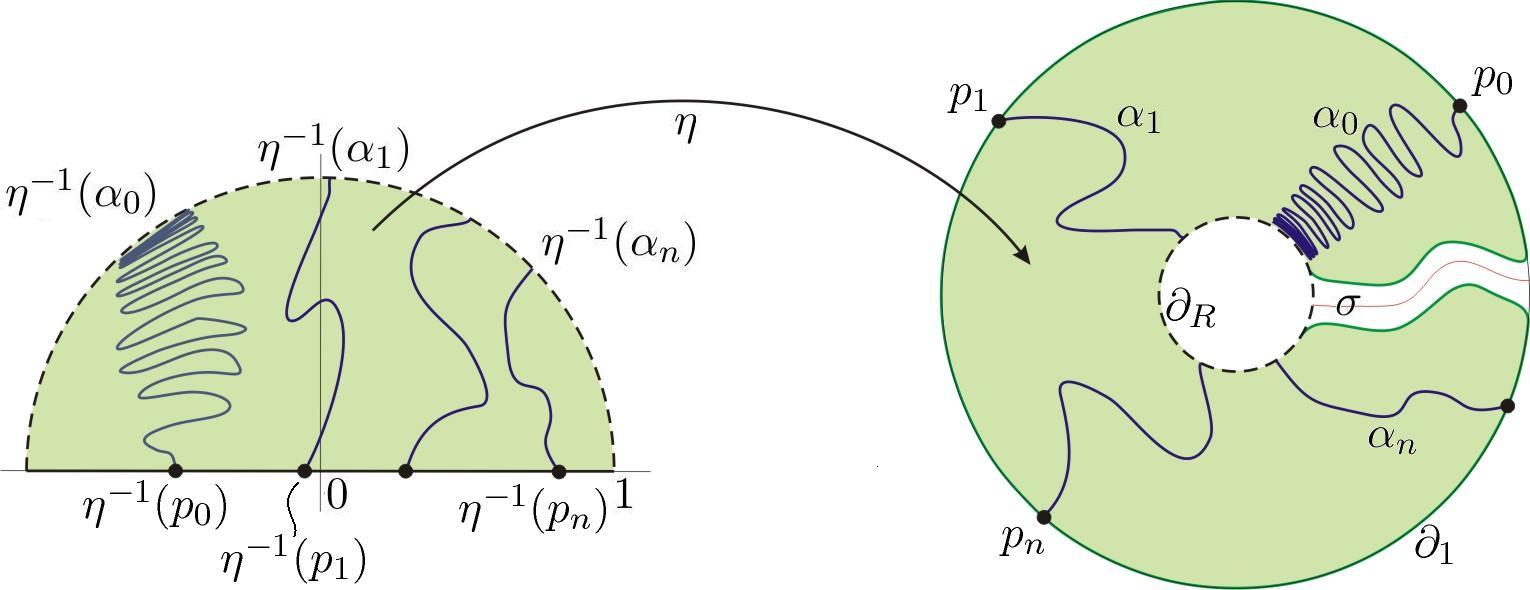}
\caption{$\eta $ conformally parameterizes the shaded region on the
right by a half disk. In fact, we will show that $\a _0$ with an
interval limit set does not occur.} \label{fig2}
\end{center}
\end{figure}
As $f$ does not have critical points in $f^{-1}(t_0)$, $F$
restricted to any of the finite number of components $\eta
^{-1}(\a_i)$ of $\widehat{f}^{-1}(t_0)$ ($0\leq i\leq n$)
monotonically parameterizes an interval on the line 
\[
L_{t_0}=\{ w=u+iv\in \C  \mid u=t_0\} . 
\]
The end points on these intervals on
$L_{t_0}$ form a finite subset, and thus, it is possible to find a
compact interval $\g $ in $L_{t_0}$ which is disjoint from the end
points in the intervals in $F(\widehat{f}^{-1}(t_0))$. Therefore,
$F^{-1}(\g)$ is compact. Hence, after replacing $A(R,1)$ by a
subend, we will assume that $F^{-1}(\g)=\mbox{\O }$. Consider the
restriction $F|_{\mbox{Int}(\Delta )}\colon \mbox{Int}(\Delta )\to
\C-\g $ to the interior of $\Delta $. $F|_{\mbox{Int}(\Delta )}$ is
essentially bounded in the sense that its image is contained in a
domain conformally equivalent to an open subset of the unit disk via
the Riemann mapping theorem. By Plessner's Theorem
(Theorem~\ref{thmplessner}), $F|_{\mbox{Int}(\Delta )}$ has angular
limits almost everywhere on $(\esf^1)^+=\{ z\in \C  \mid |z|=1,
y\geq 0\} $, and thus, $f$ has angular limits almost everywhere on
$\partial _R-\sigma (\ve )$. Clearly, by taking smaller and smaller
neighborhoods $\sigma (\ve )$, we conclude that item~(a) in
statement~2 of Theorem~\ref{th1.1} holds in this case.

We next describe the limit set $L_i$ of $\eta ^{-1}(\a_i)$,
$i=0,1,\ldots ,n$. Clearly $L_i\subset (\esf^1)^+$. Suppose that for
$i$ fixed, $\eta ^{-1}(\a _i)$ has two distinct limit points
$q_1,q_2\in (\esf^1)^+$. We claim that in this case, the subarc $I$
of $(\esf^1)^+$ whose extrema are $q_1,q_2$ is entirely contained in
$L_i$. Arguing by contradiction, suppose that there exists $s\in
\Int(I)$ which is not a limit point of $\eta ^{-1}(\a_i)$. Then
there exists a small $\de >0$ such that the radial arc $[1-\de ,1 ]s
\subset \Delta $ is disjoint from $\eta ^{-1}(\a _i)$. As $\eta
^{-1}(\a _i)$ is proper in  $\Delta $ and disjoint from $[1-\de ,1
]s$, then $\eta ^{-1}(\a _i)$ eventually lies in the one of the two
connected components, say $A$, obtained by removing $[1-\de ,1
 ]s$ from the $\de $-neighborhood of $(\esf^1)^+$ in $\Delta$.
In particular, $L_i$ fails to contain the point in $\{q_1, \, q_2\}$ which
does not lie in $\overline{A}$, which is a contradiction. This
proves our claim, and therefore, $L_i$ is either a compact subarc of
$(\esf^1)^+$ or a point.

We claim that all the limit sets $L_i$ reduce to points. Suppose on
the contrary, that some $L_i$ is a subarc of $(\esf ^1)^+$ with
nonempty interior Int$(L_i)$. Then almost everywhere on Int$(L_i)$,
the holomorphic function $F$ has angular limits which correspond to
the end point $*\in L_{t_0}\cup \{ \infty \} $ of $F(\eta
^{-1}(\a_i))$. In this case, $F$ would have the constant angular
limit $*$ on a set of positive measure of $(\esf^1)^+$, thereby
contradicting Privalov's Theorem~\cite{priv1} since $F$ is not
constant. This contradiction implies that $L_i$ is a single point in
$(\esf^1)^+$ for all $i=0,1,\ldots, n$, and so, each end of
$f^{-1}(t_0)$ has a unique limit point in $\partial_R$. Since we are
assuming $R>0$, then the harmonic measure of $\partial _R$ is
positive, which clearly implies item~(b) in statement~2
of Theorem~\ref{th1.1}. Now the proof of Assertion~\ref{a2.2} is
complete.
\end{proof}

Recall that the collection $\{ \a _0,\a _1,\ldots ,\a _n\}
=f^{-1}(t_0)$ of pairwise disjoint arcs is cyclically ordered. For
$k=0,1,\ldots, n$, let $D_k$ be the union of $\a _k\cup \a _{k+1}$
with  the component of $A(R,1)-f^{-1}(t_0)$ whose boundary contains
the arcs $\a_k, \a_{k+1}$ (in the case of $k=n$ we identify $\a _0$
with $\a _{n+1}$). Note that $\partial D_k$ is a connected open arc,
and $D_k$ is simply connected. We will call each $D_k$ a
(topological) {\it sector} of $A(R,1)$. In the next result we study
the conformal character of the sectors $D_k$.

\begin{assertion}
\label{a2.3} Suppose that $R>0$ and that the limit set
$L(f^{-1}(t_0))$ equals $\partial_R$. Then, for each $k=0,1,\ldots
,n$ the (topological) sector $D_k$ is conformally equivalent to the
closed upper halfplane in $\C $ (the conformal diffeomorphism
between $D_k$ and $\{ a+bi \mid b\geq 0 \}$ fails to be conformal
only at the starting points $p_k, \, p_{k+1}$ of $\a_k, \,
\a_{k+1}$).
\end{assertion}
\begin{proof}
The argument is again by contradiction. Assume that for some
$k=0,1,\ldots, n$ the assertion fails. By the Uniformization Theorem
and boundary regularity of holomorphic maps, then there exists an
bijective map $\phi _k\colon \Delta =\{z=x+iy \in \C \mid |z|<1, \
y\geq 0\} \to D_k$ which is a conformal immersion except at the two
points of $\partial \Delta$ that correspond to the end points $p_k,
\, p_{k+1}$ of $\a_k$ and $\a_{k+1}$ in $\partial _1$.
Since $\phi _k$ is a bounded holomorphic function,
Theorem~\ref{thmplessner} insures that we can find distinct points
$\xi _1, \xi _2\in (\esf ^1)^+$  such that $\phi _k$ has an angular
limit at $\xi _i$, $i=1,2$.  Consider a pair of smooth disjoint arcs
$\be _1,\be _2\subset \Delta $, each joining a point of
$\phi _k^{-1}[(\partial D_k\cap \partial _1)-\{ p_k,p_{k+1}\} ]\subset
\partial \Delta \cap \{ y=0\} $ to
one of the points $\xi _1,\xi _2$, and transverse to $(\esf^1)^+$.
Let $D'_k$ be the subdomain
of $\Delta $ whose boundary contains both $\be _1,\be _2$, see
Figure~\ref{fig2a}.
\begin{figure}
\begin{center}
\includegraphics[height=5cm]{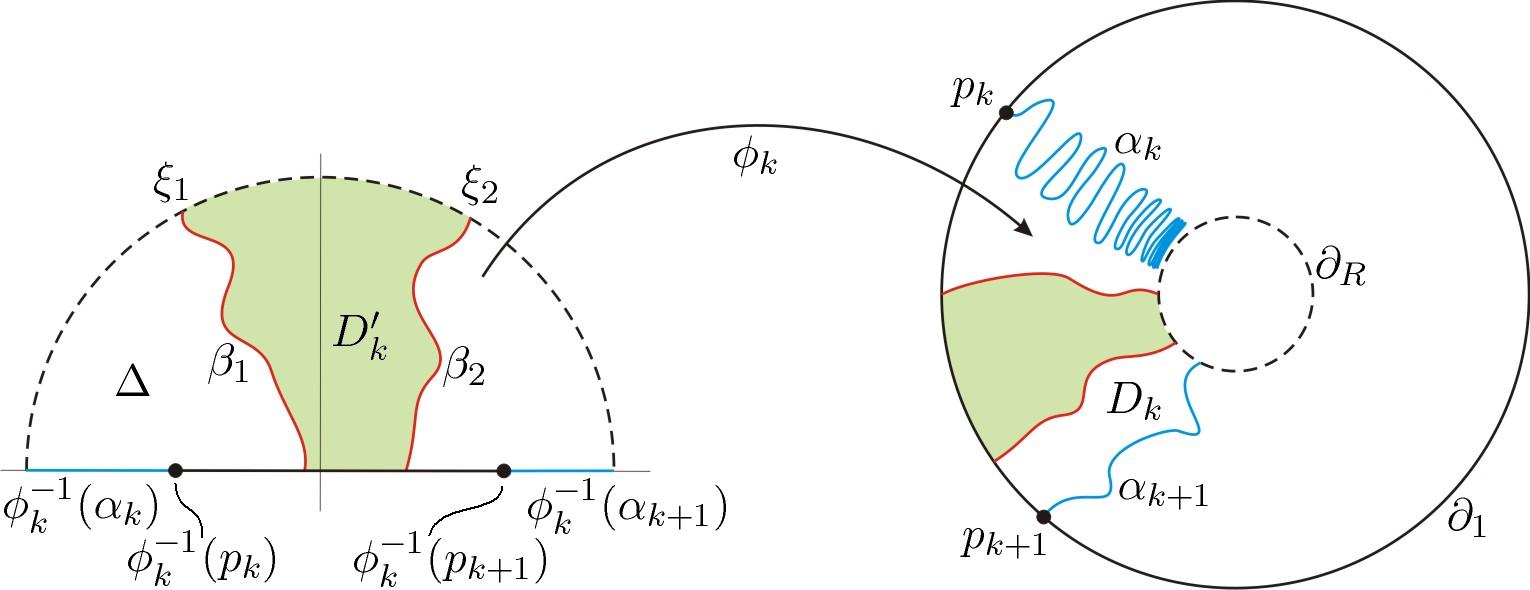}
\caption{The shaded domain on the right is bounded by curves with distinct
limiting end points.} \label{fig2a}
\end{center}
\end{figure}
Note that the (single) limit point of $\phi _k(\be _1)$ is different
from the limit point of $\phi _k(\be _2)$: otherwise for every point
$\xi '$ in the open arc $(\esf^1)^+\subset \partial \Delta $ with
extrema $\xi_1,\xi _2$, the angular limit of $\phi _k$ at $\xi '$
exists and is equal to the common limit point of $\phi _k(\be _i)$,
$i=1,2$. By Privalov's Theorem
(Theorem~\ref{thmPrivalov}), this would lead to a contradiction
since $\phi _k$ is not constant.

Also note that $\phi _k(D'_k)$ is a subdomain of $A(R,1)$ whose
boundary consists of $\phi _k(\be _1)$, $\phi _k(\be _2)$, an arc
contained in $\partial _1$ and an open arc $\de \subset \partial_R$, and
that every point of $\delta $ is a positive distance
from the domain $A(R,1)-D_k$. Therefore, $\delta $ is disjoint from
the limit set of $f^{-1}(t_0)$, which contradicts one of our
hypotheses. Now the proof of the assertion is complete.
\end{proof}

We continue with our proof of Theorem~\ref{th1.1} under the assumption  that $R>0$. Let $\Sigma$ be the flat surface
with boundary $A(R,1)-\a_0$, and let $\Sov$ denote the
simply connected flat surface with boundary obtained after attaching
to $\Sigma $ two ``copies'' of $\a _0$. Note that $\Sov$ has
connected boundary, which consists of the arc $\partial_1-\a_0$
together with the two copies of $\a_0$. For each $k\in
\{0,1,\ldots,n\}$, we will consider the "lift" $\Dhat_k$ of the
topological sector $D_k$ to $\Sov$. Let $\widehat{f}\colon \Sov \to
\R$ be the harmonic function given by the "lift" of $f$ to $\Sov$
(clearly $\widehat{f}$ takes equal values at corresponding points in
the copies of $\a _0$ in $\partial \Sov$). Since $\Sov$ is
simply connected, $\widehat{f}$ admits a well-defined harmonic
conjugate function $\widehat{f}^*$ on $\overline{\Sigma }$, and so,
the function $F=\widehat{f}+i\widehat{f}^*\colon\Sov\to \C $ is
holomorphic.

Note that each sector $\Dhat_k\subset \Sov$ has its image
$F(\Dhat_k)$ in one of the two closed halfspaces
\[
\C^+(t_0)=\{w=u+iv\in\C \mid u\geq t_0\} ,
\quad \C^{-}(t_0)=\{u+iv  \mid u\leq t_0\} .
\]
Let $L_{t_0}=\{u+iv\in \C \mid u=t_0\}$ and note that
$F^{-1}(L_{t_0})$ is the finite ordered set of curves
$\{\widehat{\a}_0, \widehat{\a}_1, \widehat{\a}_2, \ldots,
\widehat{\a}_n, \widehat{\a}'_0\}$ corresponding to the cyclically
ordered set of arcs $\{\a_0,\a_1,\ldots, \a_n\}=f^{-1}(t_0)$; also
note that $\partial \Dhat_n$ contains the arcs $\widehat{\a}_n$ and
$\widehat{\a}_0'$ in its boundary. After reindexing, we can assume
that $F(\Dhat_k)\subset \C^+(t_0)$ for $k$ even and
$F(\Dhat_k)\subset \C^-(t_0)$ for $k$ odd. For $k=0,1,\ldots ,n$,
let $w_k\in L_{t_0}\cup \{ \infty \} $ be the end point of the
half-open arc $F|_{\widehat{\a }_k}$, and let $w_{n+1}\in
L_{t_0}\cup \{ \infty \} $ be the corresponding end point of
$F|_{\widehat{\a }'_0}$.

\begin{assertion}
\label{ass3.4bis}
If $R>0$, then item~2 of Theorem~\ref{th1.1} holds.
\end{assertion}
\begin{proof}
By Assertion~\ref{a2.2}, we only need to arrive to a contradiction
provided that $L(f^{-1}(t_0))=\partial _R$. So assume
$L(f^{-1}(t_0))=\partial _R$. We first check that with the notation
above, then $w_k=w_{k+1}$ for all $k=0,1,\ldots ,n$. Consider the
restriction $F|_{\widehat{D}_k}$, which is a holomorphic function
whose image is contained in, say, $\C ^+(t_0)$. Since we are
assuming $L(f^{-1}(t_0))=\partial _R$, Assertion~\ref{a2.3} insures
that there exists a bijective map $\psi _k\colon \{ a+ib\in \C \mid
b\geq 0\} \to \widehat{D}_k$ which is conformal everywhere on the
closed upper halfplane except at the two points $q_k,q_{k+1}$ in $\{
b=0\} $ which correspond respectively to the starting points
$\widehat{p}_k,\widehat{p}_{k+1}$ of $\widehat{\a }_k,\widehat{\a
}_{k+1}$, respectively. Schwartz's reflection principle applied to
the restriction of  $F\circ \psi _k$ to $\{ a+ib \mid b\geq 0\}
-[q_k,q_{k+1}]$ (here $[q_k,q_{k+1}]$ denotes the closed interval in
the real axis with extrema $q_k,q_{k+1}$), and produces a
holomorphic function $\widetilde{F}_k \colon \C -[ q_k,q_{k+1}]\to
\C $ such that
\begin{itemize}
\item $\widetilde{F}_k$ extends continuously to the metric
completion ${\mathcal C}$ of $\C -[ q_k,q_{k+1}]$ (note that ${\mathcal C}$
is conformally $A(0,1)$). We denote  these extensions by the same
symbols $\widetilde{F}_k$.
\item $\widetilde{F}_k$ maps $\R -[
q_k,q_{k+1}]$ into the line $L_{t_0}$ and maps each of the two
copies of the interval $[ q_k,q_{k+1}]$, when considered inside the
boundary of ${\mathcal C}$, into the union of $F(\partial \widehat{D}_k
\cap \partial _1)$ and its reflected image with respect to
$L_{t_0}$. Furthermore, the preimage by $\widetilde{F}_k$ of every
point in $L_{t_0}$ consists of at most two points in the real line
$\{ b=0\} $, and for some point in $L_{t_0}$, its preimage by
$\widetilde{F}_k$ consists of at most one point.
\end{itemize}
By
Picard's theorem, $\widetilde{F}_k$ extends meromorphically across
$\infty $, and its extension is a conformal diffeomorphism from a
neighborhood $U_k$ of $\infty $ in $\C \cup \{ \infty \} $ into its
image. In particular, the limit point $w_k$ of
$\widetilde{F}_k\left( \psi _k^{-1}(\widehat{\a }_k) \right) $
equals the limit point $w_{k+1}$ of $\widetilde{F}_k\left( \psi
_k^{-1}(\widehat{\a }_{k+1}) \right) $, as desired.

We now consider the special case where $w_0=\infty $. In this case,
the pullback by $\widetilde{F}_k|_{U_k}$ of the complete flat metric
$|dw|^2$ is a complete flat metric on $U_k$. Furthermore, for each
$k$ the equality $F^*|dw|^2= (\psi _k^{-1})^*\left(
\widetilde{F}_k|_{U_k\cap \{ b\geq 0\} }\right) ^*|dw|^2$ holds on
some end representative $E_k$ of $D_k$. Therefore $F$ induces under
pullback a complete flat metric on $\cup _{k=0}^nE_k$, which is an
end representative of $\overline{\Sigma }$. Clearly this complete
flat metric on this end representative of $\overline{\Sigma }$
descends to a complete flat metric on an end representative of
$A(R,1)$ when $w_0=\infty $. This contradicts the assumption that
$R$ is positive (because any such complete flat annulus has
quadratic area growth by a direct application of the Gauss-Bonnet
formula together with the first and second variation of area formulas, and
the fact that annular ends with at most quadratic  area growth are parabolic,
see~\cite{gri1}).

On the other hand if $w_0$ is finite, then the arguments in the last
paragraph apply to the holomorphic function
$\frac{1}{\widetilde{F}_k-w_0}$ and lead to a similar contradiction.
This finishes the proof of Assertion~\ref{ass3.4bis}.
\end{proof}

In order to complete the proof of Theorem~\ref{th1.1} it remains to
demonstrate item~1 of the theorem. So from now on suppose $R=0$. As we did
just after Assertion~\ref{a2.1}, we can assume that $f\colon
A(0,1)\to \R$ is harmonic and analytic (up to and including the
boundary $\partial _1$) and $f^{-1}(t_0)$ is a cyclically ordered
finite collection $\{\a_0, \a_1, \ldots, \a_n\}$ of pairwise
disjoint, properly embedded arcs transverse to $\partial_1$, and
each arc has its starting point in $\partial_1$ and limits to $z=0$.
Also the arguments just before Assertion~\ref{a2.3} lead us to
define the topological sectors $D_k$, $k=0,1,\ldots ,n$, each one
being the union of $\a _k\cup \a _{k+1}$
with the component of
$A(0,1)-f^{-1}(t_0)$ whose boundary contains the arcs $\a_k,
\a_{k+1}$ (with $\a _0=\a _{n+1}$ if $k=n$), and with an arc in
$\partial _1$. Note that these sectors $D_k$ are still parabolic in
our new setting of $R=0$, since the conformal structure of the
annulus $A(0,1)$ is parabolic and the sectors $D_k$ are then proper
subdomains of a parabolic surface.

Repeating the arguments before Assertion~\ref{ass3.4bis}, we cut
$A(0,1)$ along $\a _0$ and then reattach the cutting curve twice to
obtain a simply connected surface $\overline{\Sigma }$, a
holomorphic function $F=\widehat{f}+i\widehat{f}^*\colon
\overline{\Sigma }\to \C $ and a finite, ordered set of arcs
$\{\widehat{\a}_0, \widehat{\a}_1, \widehat{\a}_2, \ldots,
\widehat{\a}_n, \widehat{\a}'_0\}$ which correspond to
$\{\a_0,\a_1,\ldots, \a_n\}=f^{-1}(t_0)\subset A(0,1)$. By the
arguments in the proof of Assertion~\ref{ass3.4bis}, the holomorphic
differential $dF$ on $\overline{\Sigma }$ descends to the
holomorphic differential $\omega =df+idf^*$ on $A(0,1)$, which
extends to a meromorphic differential $\widetilde{\omega }$ on
$\overline{\D }=\overline{A}(0,1)$.

In order to prove item~1(a), note that if $\widetilde{\omega }$
is holomorphic at $0\in \overline{A}(0,1)$, then $d\widehat{F}$ is
holomorphic at $\infty $ (with the same notation as in the proof of
Assertion~\ref{ass3.4bis}). This implies that $F$ is holomorphic at
$0$ and so, $f$ is bounded on $E$. Reciprocally, if $f$ is bounded
then $F$ is bounded as well, which implies $\widetilde{\omega }=dF$
is holomorphic.

Finally we prove item~1(b) of Theorem~\ref{th1.1}. Since
length$(\a )=\int _{\a }|*df|=\int _{\a }|df|<\infty $, then the
common limit point $w_0=\ldots =w_n=w_{n+1}$ corresponds a finite
point. In this case, the arguments in the last paragraph demonstrate
that $f$ is bounded.

\begin{remark}
 {\rm
Suppose $R>0$ and $f\colon A(R,1)\to \R$ is a nonconstant harmonic
function with angular limits almost everywhere on $\partial_R$. Then
the arguments in the proof of Theorem~\ref{th1.1} can also be
applied to demonstrate the following result: For every $t\in \R$,
each proper arc (piecewise smooth) in the boundary of a component of
$\{z\in A(1,R) \mid f(z)\geq t\}$ or of $\{z\in A(1,R) \mid f(z)\leq
t\}$ has a unique limit point in $\partial_R$. In particular, each
nonlimit end of the $1$-complex $f^{-1}(t)$ has a unique limit point
in $\partial_R$. }
\end{remark}

\section{The proof of Corollary~\ref{c1.2}.}

In this section, we will apply the following theorem of Collin,
Kusner, Meeks and Rosenberg~\cite{ckmr1} to prove
Corollary~\ref{c1.2}.

\begin{theorem}
\label{t3.1} If $X\colon \Sigma \to \rth$ is a properly immersed
minimal surface with  boundary (possibly empty) which is contained
in a halfspace of $\rth$, then $\Sigma$ is parabolic.
\end{theorem}

Let $R\in [0,1)$. Suppose $X\colon A(R,1)\to \R^3$ is a conformal,
proper minimal immersion such that $X^{-1}(P)$ has a finite number
of ends for some horizontal plane $P\subset \R^3$ at height $t_0\in
\R $.

We claim that $R=0$. Otherwise $R>0$, $A(R,1)$ is an annular end of
finite type for $x_3$ and some component of $x_3^{-1}(( -\infty
,t_0])$ or $x_3^{-1}([t_0,\infty ))$ is hyperbolic by item~2 of
Theorem~\ref{th1.1}. But such a component must be parabolic by
Theorem~\ref{t3.1}, since $X$ restricted to this component is a
properly immersed minimal surface contained in a halfspace of
$\rth$. Hence $R=0$.

By item~1 of Theorem~\ref{th1.1}, the holomorphic one-form
$dx_3+idx_3^*$ extends to a meromorphic one-form on $\overline{\D
}=\overline{A}(0,1)$. By regularity of the induced metric, the
meromorphic Gauss map $g\colon A(0,1)\to \C \cup \{ \infty \} $ of
$X$ achieves the values  $0,\infty $, corresponding to the north and
south poles of $\esf^2$ and equal to the unit normals of $P$,  only
a finite number of times. Hence, $g$ misses $0,\infty $ on $A(0,R')$
for some $R'\in (0,1]$. For
some $k\in \Z $, $z^{-k}g$ induces the zero map from $\pi_1(A(0,R'))$ to
$\pi_1(\C - \{0\})$  and thus, by elementary covering space theory, $z^{-k}g(z)=e^{H(z)}$
for some holomorphic function $H$ on $A(0,R')$. This completes the
proof of the main statement of Corollary~\ref{c1.2}.

We next prove the equivalence between items 1--3 of
Corollary~\ref{c1.2}. The only implication which is not obvious is
that 3 implies 1, so assume 3 holds. The main statement
of Corollary~\ref{c1.2} applied to each of the planes $P_1, P_2$
implies that on some end $A_{R_1}=\{ z\in A(0,1)\mid 0<|z|\leq
R_1<1\} $ of $A(0,1)$, the Gauss map $g$ misses the four values
corresponding to the two pairs of antipodal points of $\esf^2$ which
are orthogonal to $P_1$ or $P_2$. Since $A_{R_1}$ is conformally a
punctured disk, Picard's theorem can be applied to $g$ and gives
that $g$ extends across the puncture to a meromorphic function on
$\overline{\D}=\overline{A(0,1)}$. Hence, $X$ has finite total
curvature.

Finally suppose that $\int _{\a }\frac{\partial x_3}{\partial \eta
}\, ds$ is finite for some end representative of $X_3^{-1}(P)$. Then
item~1.(b) of Theorem~\ref{th1.1} implies that $x_3$ is bounded
on the minimal end $E=X(A(0,1))$. In particular, $E$ is contained in
a horizontal slab. In this setting, Lemma~2.2 in~\cite{ckmr1}
insures that $E$ has quadratic area growth. Since the Gaussian
curvature of $E$ is nonpositive, a standard application of the first
and second variation formulas for area imply that $E$ has finite
total curvature. In this situation, it is well-known that $E$ is
asymptotic to $P$ with finite multiplicity. Now Corollary~\ref{c1.2}
is proved.
\newpage

\center{William H. Meeks, III at profmeeks@gmail.com\\
Mathematics Department, University of Massachusetts, Amherst, MA 01003}
\center{Joaqu\'\i n P\'{e}rez at jperez@ugr.es\\
Department of Geometry and Topology and Institute of Mathematics
(IEMath-GR), University of Granada, 18071, Granada, Spain}

\bibliographystyle{plain}
\bibliography{bill}

\begin{thebibliography}{1}

\bibitem{allo1}
A.~Alarc\'on and F.~J. L\'opez.
\newblock On harmonic quasiconformal immersions in $\mathbb{R}^3$.
\newblock {\em Transactions of the AMS}, 365(4):1711--1742, 2013.
\newblock MR3009644, Zbl 1320.53076.

\bibitem{ckmr1}
P.~Collin, R.~Kusner, W.~H. Meeks~III, and H.~Rosenberg.
\newblock The geometry, conformal structure and topology of minimal surfaces
  with infinite topology.
\newblock {\em J. Differential Geom.}, 67:377--393, 2004.
\newblock MR2153082, Zbl 1098.53006.

\bibitem{gri1}
A.~Grigor{'}yan.
\newblock Analytic and geometric background of recurrence and non-explosion of
  {B}rownian motion on {R}iemannian manifolds.
\newblock {\em Bull. of {A}.{M}.{S}}, 36(2):135--249, 1999.
\newblock MR1659871, Zbl 0927.58019.

\bibitem{lop4}
F.~J. L{\'o}pez.
\newblock Some {P}icard theorems for minimal surfaces.
\newblock {\em Trans. Amer. Math. Soc.}, 356(2):703--733, 2004.
\newblock MR2022717 (2004j:53017), Zbl 1121.53007.

\bibitem{mpe3}
W.~H. Meeks~III and J.~P\'{e}rez.
\newblock Embedded minimal surfaces of finite topology.
\newblock Preprint at http://arxiv.org/pdf/1506.07793v1.pdf.

\bibitem{mr7}
W.~H. Meeks~III and H.~Rosenberg.
\newblock Maximum principles at infinity.
\newblock {\em J. Differential Geom.}, 79(1):141--165, 2008.
\newblock MR2401421, Zbl pre05285650.

\bibitem{ple1}
A.~I. Plessner.
\newblock \"{U}ber das verhalten analytischer functionen am rande ihres
  definitionbereichs.
\newblock {\em J. Reine Angew. Math.}, 158:219--227, 1927.
\newblock JFM 53.0284.01.

\bibitem{priv1}
I.~I. Privalov.
\newblock {\em Randeigenschaften analytischer Funktionen. 2. \"Uberarb. und
  erg. Aufl.}
\newblock Hochschulb\"ucher f\"ur Mathematik. 25. Berlin: VEB Deutscher Verlag
  der Wissenschaften VIII, 247 S. mit 10 Abb., 1956.
\newblock Zbl 0073.06501.

\end{thebibliography}
\end{document}